\documentclass[12pt]{article}
\usepackage{amsmath}
\usepackage{amsfonts}
\usepackage{amssymb}
\usepackage[dvips,unicode]{hyperref}
\usepackage{graphicx}
\usepackage{cite}

\begin{document}

\title{ Construction of the surfaces with horizontal tangential planes at given
points }
\author{I.P. Smirnov and V.G. Burdukovskaya\\
Institute of Applied Physics RAS, \\
46 Ul'yanova Street, Nizhny Novgorod, Russia}
\date{}
\maketitle

\begin{abstract}
The problem of construction of the surfaces with given sets of the points with
horizontal tangential planes is considered. Such considerations are of
interest in the problem of computer simulations of the waved ocean surfaces.

\end{abstract}

\section{Statement of the problem}

\subsection{Physical statement}

Various problems of optical/acoustical monitoring of wavy water surface are
related with simulations of the surfaces with given properties \cite{l2,l4,l5}%
. The example of it is the problem of construction of two-dimension surface
having horizontal target planes at the given points (Luchinin A.G., 2003).
Variants of solutions of this problem are presented in this article.

\subsection{Mathematical formulation}

\textbf{The main problem A. }\emph{Let }$\mathfrak{M=}\left\{  \mathbf{x}%
_{i},\ i=\overline{0,N}\right\}  $\emph{ be the finite subset of the given
open set }$\Omega\subset R^{n}\emph{,}n\geq1.$\emph{ Construct  a function
}$F\left(  \mathbf{x}\right)  $\emph{ of the given class }$\mathfrak{F}\in
C^{1}\left(  \Omega\right)  $ \emph{possessing zero gradient }$\nabla F\left(
\mathbf{x}\right)  =\mathbf{0}$ \emph{at all points of }$\mathfrak{M}$\emph{
.}

We study also the next \textit{specifications} of the main problem:

\begin{description}
\item[B.] \emph{the values of the function at the points }$\mathfrak{M}$\emph{
are given (}$F\left(  \mathbf{x}_{i}\right)  =z_{i},\ i=\overline{0,N}%
$\emph{);}

\item[C.] \emph{gradient of the function takes zero values only in
}$\mathfrak{M}$\emph{.}
\end{description}

It is clear that sets of the solutions of \textbf{A}-\textbf{C} are infinite.
In general, solutions of the main problem \textbf{A} do not solve problems
\textbf{B,C}. The main goal of the next research is to determine the existence
of the solution of a problem in a given class of functions $\mathfrak{F}$.

\section{Analytical approach}

\label{anal}

Let $\mathfrak{F}=C^{k}\left(  \Omega\right)  $, $1\leq k\leq\infty$. For this
class of functions we can easily construct an analytic solution of the main
problem\textbf{. }One of possible constructions is described below. Let
$\left\vert \cdot\right\vert $ be Euclidean norm in $R^{n}$. Consider for
every point $\mathbf{x}_{i}$ $\in\mathfrak{M}$ the open sphere%
\[
S_{i}\equiv\left\{  \mathbf{x}:\left\vert \mathbf{x-x}_{i}\right\vert
<r_{i}\right\}  \subset\Omega
\]
with the radius $r_{i}>0$ and center point $\mathbf{x}_{i}$. Let $a_{i}>1$,
$c_{i}$ be arbitrary values, $\varphi_{i}\left(  \cdot\right)  \in
C^{k}\left(  R^{n}\right)  $ be arbitrary functions such that $\varphi
_{i}\left(  \mathbf{0}\right)  =1$, $\nabla\varphi_{i}\left(  \mathbf{0}%
\right)  =0$. Let%
\begin{gather}
\chi_{i}\left(  \mathbf{x}\right)  \equiv\left\{
\begin{array}
[c]{ll}%
a_{i}^{\frac{\left\vert \mathbf{x}\right\vert ^{2}}{\left\vert \mathbf{x}%
\right\vert ^{2}-r_{i}^{2}}}, & \left\vert \mathbf{x}\right\vert <r_{i},\\
0, & \left\vert \mathbf{x}\right\vert \geq r_{i},
\end{array}
\right.  \nonumber\\
F\left(  \mathbf{x}\right)  =%
{\displaystyle\sum\limits_{i=0}^{N}}
c_{i}\chi_{i}\left(  \mathbf{x-x}_{i}\right)  \varphi_{i}\left(
\mathbf{x-x}_{i}\right)  .\label{F_sol}%
\end{gather}

It easy to prove the following statements.

\textbf{Proposition 1.} \emph{For every set of sufficiently small values
}$r_{i}$\emph{ (such that }$\mathbf{x}_{i}\notin$\emph{ }$S_{j}$, $j\neq
i$\emph{) and for any values }$c_{i}$\emph{ the function} \emph{(\ref{F_sol})
is a solution of the problem \textbf{A}. In addition, if}

\begin{enumerate}
\item $c_{i}=z_{i},\ i=\overline{0,N}$\emph{, then (\ref{F_sol}) gives a
solution of the problem\textbf{ B};}

\item \emph{all spheres }$S_{j}$\emph{ do not intersect, then (\ref{F_sol})
gives a solution of the problem\textbf{ C}.}
\end{enumerate}

\textbf{Proposition 2. }\emph{Let }$n=1$\emph{, }$x_{0}<x_{1}<\ldots<x_{N}%
$\emph{ . Then under the conditions}%
\[
r_{i}\leq\left\vert x_{i}-x_{i-1}\right\vert ,~c_{i-1}c_{i}<0,\ i=\overline
{1,N},
\]
\emph{the function (\ref{F_sol}) gives a solution of the problem \textbf{C}.}

The peculiarity of the described construction is the provided by multipliers
$\chi_{i}\left(  \mathbf{x}\right)  $ finiteness of the terms of the sum
(\ref{F_sol}); this provides zero values of gradient at the points of
$\mathfrak{M}$. Note that using the corresponding values of gradients
$\nabla\varphi_{i}\left(  \mathbf{0}\right)  $ we can assume these gradients
not only zero but also any given, so we can use the construction (\ref{F_sol})
for solving of more complicated problems to construct of surfaces with given
slopes at the given points.

\begin{figure}[ptb]%
\centering
\includegraphics[
height=7.8233cm,
width=10.3941cm
]%
{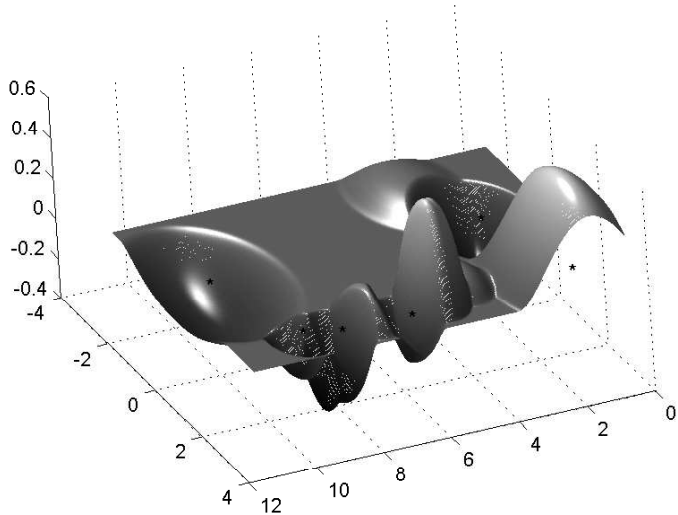}%
\caption{The surface (\ref{F_sol}) for the set of $10$ random points marked on
the graph by stars.}%
\label{analyt}%
\end{figure}

From physical point of view the solutions (\ref{F_sol}) can simulate only
oscillation of too viscous liquid. An example of the numerical simulation by
means of formula (\ref{F_sol}) for $a_{i}=e$, $\varphi_{i}\left(
\mathbf{x}\right)  \equiv1$, $10$ random points $\mathfrak{M}$ and random
coefficients $c_{i}$, is presented on Fig. \ref{analyt}.%

\section{Spline approach}

\label{spline}

In this section we use splines of various types to solve problems
\textbf{A}-\textbf{C} for $n=1,2$. Polynomial splines of the 4-th order (see
sec. \ref{spline14}, \ref{spline24}) make it possible to solve problems
\textbf{A} and \textbf{B}. As for the problem \textbf{C} then in the case
$n=2$ spline approach is useless if only the set $\mathfrak{M}$ does not have
the form of the product
\begin{equation}
\mathfrak{M}=X\times Y,\ \ \left\{
\begin{array}
[c]{l}%
X=\left\{  x_{k}\right\}  ,\ k=\overline{0,K},\\
Y=\left\{  y_{l}\right\}  ,\ \ l=\overline{0,L}.
\end{array}
\right.  \label{Matr}%
\end{equation}
But in the case $n=1$ the solution of the problem \textbf{C} can be found in
the class of trigonometric splines (sec. \ref{spline1t}).

\subsection{Polynomial splines}

\subsubsection{One-dimensional problem ($n=1$)}

\label{spline14}

Consider the following problem \textbf{B}: \emph{for given sets }%
$a=x_{0}<x_{1}<\ldots<x_{N}=b$\emph{\ and }$\left\{  z_{i},i=\overline
{0,N},\right\}  $\emph{ build a spline of the 4-th order }$F\left(
\cdot\right)  \in C^{1}\left[  a,b\right]  $\emph{\ such that}%
\begin{equation}
\left\{
\begin{array}
[c]{l}%
F\left(  x\right)  =s_{i}\left(  x\right)  ,\ x\in\left[  x_{i-1}%
,x_{i}\right]  ,\ i=\overline{1,N},\\
s_{i}\left(  x\right)  \equiv a_{i}+b_{i}\left(  x-x_{i}\right)  +c_{i}\left(
x-x_{i}\right)  ^{2}/2!+d_{i}\left(  x-x_{i}\right)  ^{3}/3!+e_{i}\left(
x-x_{i}\right)  ^{4}/4!,\\
F\left(  x_{i-1}\right)  =z_{i-1},\ \ \ F\left(  x_{i}\right)  =z_{i},\\
F^{\prime}\left(  x_{i-1}\right)  =0,\ \ F^{\prime}\left(  x_{i}\right)  =0.
\end{array}
\right.  \label{s1}%
\end{equation}

\begin{figure}[tbh]%
\centering
\includegraphics[
trim=0.000000cm 0.000000cm 0.000000cm -0.468557cm,
height=7.7104cm,
width=9.637cm
]%
{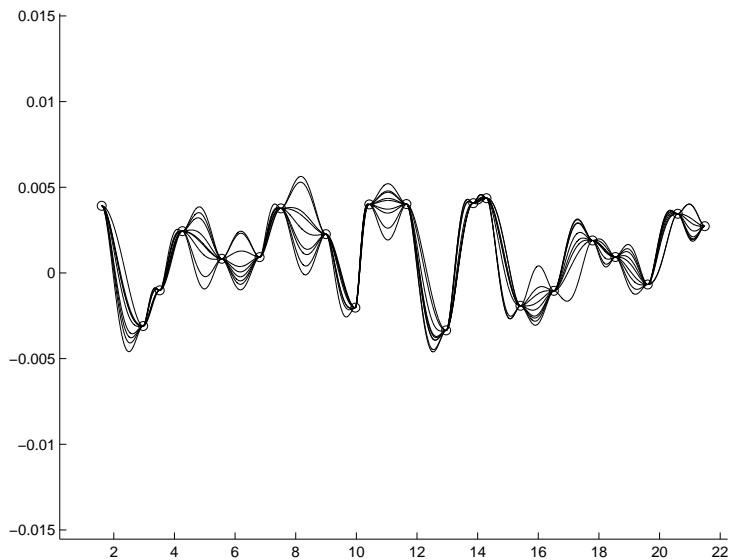}%
\caption{Examples of smooth polynomial splines of the 4-th order solving the
problem \textbf{B}. The choice rule for the coefficients : $c_{i}=\xi$, where
$\xi$ is a random value  uniformly distributed in the interval $\left[
-0.05,0.05\right]  $. Points $\left(  x_{i},z_{i}\right)  $ are marked on the
graphs by circles.}%
\label{fig2}%
\end{figure}

Now we prove that such splines exist. We have for the coefficients of the
spline%
\[
\left\{
\begin{array}
[c]{l}%
F\left(  x_{i}\right)  =a_{i}=z_{i},\\
F^{\prime}\left(  x_{i}\right)  =b_{i}=0,\\
F\left(  x_{i-1}\right)  =a_{i}+b_{i}\Delta x_{i}+c_{i}\Delta x_{i}%
^{2}/2+d_{i}\Delta x_{i}^{3}/6+e_{i}\Delta x_{i}^{4}/24=z_{i-1},\\
F^{\prime}\left(  x_{i-1}\right)  =b_{i}+c_{i}\Delta x_{i}+d_{i}\Delta
x_{i}^{2}/2+e_{i}\Delta x_{i}^{3}/6=0,\\
\Delta x_{i}\equiv x_{i-1}-x_{i},\ i=\overline{1,N}.
\end{array}
\right.
\]
So we obtain the equations%
\[
\left\{
\begin{array}
[c]{l}%
c_{i}=-d_{i}\Delta x_{i}/2-e_{i}\Delta x_{i}^{2}/6,\\
-d_{i}\Delta x_{i}^{3}/12-e_{i}\Delta x_{i}^{4}/24=\Delta z_{i}\equiv
z_{i-1}-z_{i},\\
d_{i}+e_{i}\Delta x_{i}/3=-2c_{i}/\Delta x_{i},\\
d_{i}+e_{i}\Delta x_{i}/2=-12\Delta z_{i}/\Delta x_{i}^{3},
\end{array}
\right.
\]
and we can evaluate all coefficient by $c_{i}$,
\begin{equation}
\left\{
\begin{array}
[c]{l}%
a_{i}=z_{i},\\
b_{i}=0,\\
c_{i}=c_{i},\\
d_{i}=-6\left(  c_{i}-4\Delta z_{i}\Delta x_{i}^{-2}\right)  /\Delta x_{i},\\
e_{i}=12\left(  c_{i}-6\Delta z_{i}\Delta x_{i}^{-2}\right)  /\Delta x_{i}%
^{2}.
\end{array}
\right.  \label{abcde}%
\end{equation}
Spline and its derivative are continuous both. In general, for arbitrary
coefficients $c_{i}$ we obtain smooth splines with discontinuous second
derivative . The examples of such splines are presented on Fig. \ref{fig2}%

\begin{figure}[tbh]%
\centering
\includegraphics[
trim=0.000000cm 0.000000cm 0.000000cm -0.593378cm,
height=7.7104cm,
width=9.6348cm
]%
{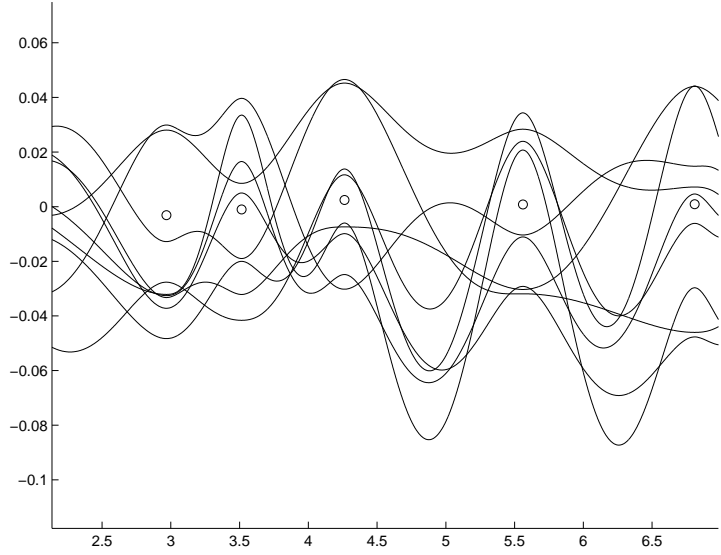}%
\caption{Examples of polynomial splines of the fourth order from
~$C^{2}\left[  a,b\right]  $ solving the problem \textbf{A.}}%
\label{fig4}%
\end{figure}

If we do not fix parameters $z_{i}$ in (\ref{abcde}) then splines solve only
the problem \textbf{A (}see Fig . \ref{fig4})%

\subparagraph{$C^{2}$-spline}

Let us now construct a spline with the continuous second derivative. For
realization of it we must require%
\[
s_{i-1}^{\prime\prime}\left(  x_{i-1}\right)  =s_{i}^{\prime\prime}\left(
x_{i-1}\right)  ,\ i=\overline{1,N},
\]
that tends to the equations%
\[
c_{i-1}=c_{i}+d_{i}\Delta x_{i}+\frac{e_{i}}{2}\Delta x_{i}^{2}=c_{i}-12\Delta
z_{i}\Delta x_{i}^{-2}.
\]
This provides the recurrence%
\begin{equation}
c_{i}=c_{i-1}+12\Delta z_{i}\Delta x_{i}^{-2},\label{rec}%
\end{equation}
so%
\[
c_{i}=c_{0}+12%
{\displaystyle\sum\limits_{j=1}^{i}}
\Delta z_{j}\Delta x_{j}^{-2},\ i=\overline{1,N}.
\]

\begin{figure}[ptb]%
\centering
\includegraphics[
height=8.6954cm,
width=11.5505cm
]%
{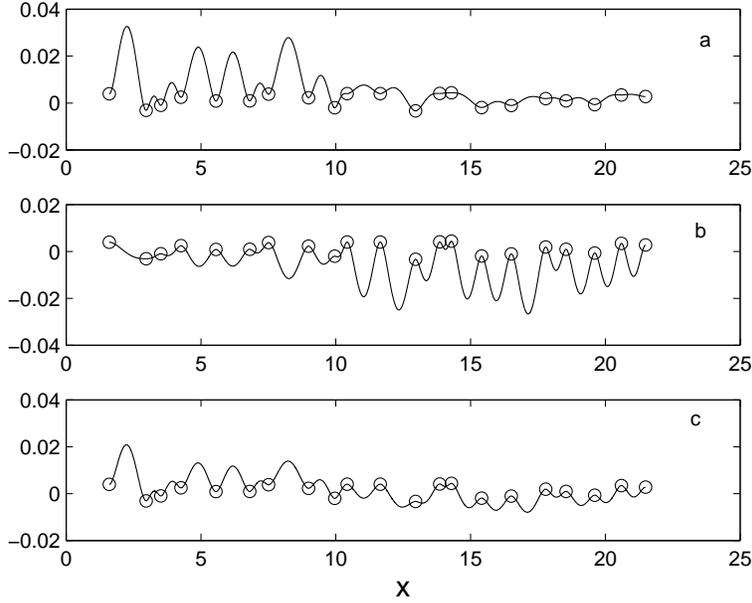}%
\caption{Polynomial splines of the forth order with additional properties: a)
all points of $\mathfrak{M}$ are points of local minima, b) all points of
$\mathfrak{M}$ are points of local maxima, c) the mean value of the second
derivative for over all points of $\mathfrak{M}$ takes its minimum.}%
\label{figabc}%
\end{figure}

Hence we obtain a one-parameter aggregate of splines of the class
$C^{2}\left[  a,b\right]  $. The parameter of it is $c_{0}$, the second
derivative of the spline at the initial point $a$. By choosing of $c_{0}$, we
can construct splines with different properties. For example, when%
\[
c_{0}=-12\min_{i}%
{\displaystyle\sum\limits_{j=1}^{i}}
\Delta z_{j}\Delta x_{j}^{-2},
\]
all points of $\mathfrak{M}$ are points of local minima (see Fig.
\ref{figabc}, a), for%
\[
c_{0}=-12\max_{i}%
{\displaystyle\sum\limits_{j=1}^{i}}
\Delta z_{j}\Delta x_{j}^{-2},
\]
all points $\mathfrak{M}$ are points of local maxima (see Fig. \ref{figabc},
b), and for%
\[
c_{0}=-\frac{12}{n}%
{\displaystyle\sum\limits_{j=1}^{i}}
\Delta z_{j}\Delta x_{j}^{-2}.
\]
the mean value of the second derivative for over all points of $\mathfrak{M}$
takes its minimum (see Fig. \ref{figabc}, c)

\subsubsection{Two-dimensional problem ($n=2$)}

\label{spline24}

Consider the following two-dimensional problem \textbf{B: } \emph{for given
set of points }$\left(  x_{k},y_{l},z_{kl}\right)  $\emph{, }$a=x_{0}%
<x_{1}<\ldots<x_{K}=b$, $c=y_{0}<y_{1}<\ldots<y_{L}=d$\emph{\ build a spline
of the 4-th order }$F\left(  \cdot\right)  \in C^{1}\left(  \left[
a,b\right]  \times\left[  c,d\right]  \right)  $\emph{\ such that}%
\begin{equation}
\left\{
\begin{array}
[c]{l}%
F\left(  x,y\right)  =s_{kl}\left(  x,y\right)  ,\ \left(  x,y\right)
\in\mathbf{P}_{kl}\equiv\left[  x_{k-1},x_{k}\right]  \times\left[
y_{l-1},y_{l}\right]  \ k=\overline{1,K},~l=\overline{1,L},\\
s_{kl}\left(  x,y\right)  =%
{\displaystyle\sum\limits_{i,j=0}^{4}}
a_{ij}^{\left(  kl\right)  }\left(  x-x_{k}\right)  ^{i}\left(  y-y_{l}%
\right)  ^{j},\\
F\left(  x_{k},y_{l}\right)  =z_{kl},\\
\triangledown F\left(  x_{k},y_{l}\right)  =0.
\end{array}
\right.  \label{s2}%
\end{equation}

By decomposition of the function into one-dimensional splines (\ref{s1}) onto
grid lines $\left(  x_{k},y_{l}\right)  ,k=\overline{1,K},~l=\overline{1,L}$,
we receive the following equations for the coefficients of the two-dimensional
spline,%
\begin{equation}%
\begin{array}
[c]{ll}%
{\displaystyle\sum\limits_{i=0}^{4}}
\left(
{\displaystyle\sum\limits_{j=0}^{4}}
a_{ij}^{\left(  kl\right)  }\Delta y_{l}^{j}\right)  \left(  x-x_{k}\right)
^{i}=%
{\displaystyle\sum\limits_{i=0}^{4}}
\gamma_{i}^{\left(  k\ l-1\right)  }\left(  x-x_{k}\right)  ^{i}, &
y=y_{l-1},\\%
{\displaystyle\sum\limits_{i=0}^{4}}
a_{i0}^{\left(  kl\right)  }\left(  x-x_{k}\right)  ^{i}=%
{\displaystyle\sum\limits_{i=0}^{4}}
\gamma_{i}^{\left(  kl\right)  }\left(  x-x_{k}\right)  ^{i}, & y=y_{l},\\%
{\displaystyle\sum\limits_{j=0}^{4}}
\left(
{\displaystyle\sum\limits_{i=0}^{4}}
a_{ij}^{\left(  kl\right)  }\Delta x_{k}^{i}\right)  \left(  y-y_{l}\right)
^{j}=%
{\displaystyle\sum\limits_{j=0}^{4}}
\beta_{j}^{\left(  k-1\ l\right)  }\left(  y-y_{l}\right)  ^{j}, &
x=x_{k-1},\\%
{\displaystyle\sum\limits_{j=0}^{4}}
a_{0j}^{\left(  kl\right)  }\left(  y-y_{l}\right)  ^{j}=%
{\displaystyle\sum\limits_{j=0}^{4}}
\beta_{j}^{\left(  kl\right)  }\left(  y-y_{l}\right)  ^{j}, & x=x_{k},
\end{array}
\label{urav2}%
\end{equation}
where we know values $\gamma_{i}^{\left(  kl\right)  },\beta_{j}^{\left(
kl\right)  }$ from formulae (\ref{abcde}); so we know the first column and the
first row of the matrix $A^{\left(  kl\right)  }\equiv\left\{  a_{ij}^{\left(
kl\right)  }\right\}  _{i,j=0}^{4}$ . Last $16$ coefficients%
\[
A^{\left(  kl\right)  }=\left\Vert
\begin{array}
[c]{ccccc}%
z_{kl} & 0 & \beta_{2}^{\left(  kl\right)  } & \beta_{3}^{\left(  kl\right)  }
& \beta_{4}^{\left(  kl\right)  }\\
0 & a_{11}^{\left(  kl\right)  } & a_{12}^{\left(  kl\right)  } &
a_{13}^{\left(  kl\right)  } & a_{14}^{\left(  kl\right)  }\\
\gamma_{2}^{\left(  kl\right)  } & a_{21}^{\left(  kl\right)  } &
a_{22}^{\left(  kl\right)  } & a_{23}^{\left(  kl\right)  } & a_{24}^{\left(
kl\right)  }\\
\gamma_{3}^{\left(  kl\right)  } & a_{31}^{\left(  kl\right)  } &
a_{32}^{\left(  kl\right)  } & a_{33}^{\left(  kl\right)  } & a_{34}^{\left(
kl\right)  }\\
\gamma_{4}^{\left(  kl\right)  } & a_{41}^{\left(  kl\right)  } &
a_{42}^{\left(  kl\right)  } & a_{43}^{\left(  kl\right)  } & a_{44}^{\left(
kl\right)  }%
\end{array}
\right\Vert
\]
are unknown. For them we have from (\ref{urav2}) the equations%
\begin{equation}
\left\{
\begin{array}
[c]{c}%
{\displaystyle\sum\limits_{j=1}^{4}}
a_{ij}^{\left(  kl\right)  }\Delta y_{l}^{j}=\gamma_{i}^{\left(
k\ l-1\right)  }-\gamma_{i}^{\left(  kl\right)  }\equiv\Delta\gamma
_{i}^{\left(  kl\right)  },\ i=\overline{1,4},\\%
{\displaystyle\sum\limits_{i=1}^{4}}
a_{ij}^{\left(  kl\right)  }\Delta x_{k}^{i}=\beta_{j}^{\left(  k-1\ l\right)
}-\beta_{j}^{\left(  kl\right)  }\equiv\Delta\beta_{j}^{\left(  kl\right)
},\ j=\overline{1,4}.
\end{array}
\right.  \label{sys8}%
\end{equation}
The matrix of this system is%
\[
D=\left\Vert
\begin{array}
[c]{cccccccc}%
\Delta y_{l} & 0 & 0 & 0 & \Delta x_{k} & 0 & 0 & 0\\
\Delta y_{l}^{2} & 0 & 0 & 0 & 0 & \Delta x_{k} & 0 & 0\\
\Delta y^{3} & 0 & 0 & 0 & 0 & 0 & \Delta x_{k} & 0\\
\Delta y_{l}^{4} & 0 & 0 & 0 & 0 & 0 & 0 & \Delta x_{k}\\
0 & \Delta y_{l} & 0 & 0 & \Delta x_{k}^{2} & 0 & 0 & 0\\
0 & \Delta y_{l}^{2} & 0 & 0 & 0 & \Delta x_{k}^{2} & 0 & 0\\
0 & \Delta y_{l}^{3} & 0 & 0 & 0 & 0 & \Delta x_{k}^{2} & 0\\
0 & \Delta y_{l}^{4} & 0 & 0 & 0 & 0 & 0 & \Delta x_{k}^{2}\\
0 & 0 & \Delta y_{l} & 0 & \Delta x_{k}^{3} & 0 & 0 & 0\\
0 & 0 & \Delta y_{l}^{2} & 0 & 0 & \Delta x_{k}^{3} & 0 & 0\\
0 & 0 & \Delta y_{l}^{3} & 0 & 0 & 0 & \Delta x_{k}^{3} & 0\\
0 & 0 & \Delta y_{l}^{4} & 0 & 0 & 0 & 0 & \Delta x_{k}^{3}\\
0 & 0 & 0 & \Delta y_{l} & \Delta x_{k}^{4} & 0 & 0 & 0\\
0 & 0 & 0 & y^{2} & 0 & x^{4} & 0 & 0\\
0 & 0 & 0 & y^{3} & 0 & 0 & x^{4} & 0\\
0 & 0 & 0 & y^{4} & 0 & 0 & 0 & x^{4}%
\end{array}
\right\Vert ,
\]
$\allowbreak$ and its rank is equal $7$. If we choose
\[
a_{11}^{\left(  kl\right)  },a_{12}^{\left(  kl\right)  },a_{13}^{\left(
kl\right)  },a_{14}^{\left(  kl\right)  },a_{21}^{\left(  kl\right)  }%
,a_{31}^{\left(  kl\right)  },a_{41}^{\left(  kl\right)  },
\]
as known values (they correspond to nonzero minor of the order $7$) and choose
other $9$ elements as free, then we obtain from (\ref{sys8}) the system%
\begin{equation}
\left\{
\begin{array}
[c]{l}%
a_{11}^{\left(  kl\right)  }\Delta y_{l}+a_{12}^{\left(  kl\right)  }\Delta
y_{l}^{2}+a_{13}^{\left(  kl\right)  }\Delta y_{l}^{3}+a_{14}^{\left(
kl\right)  }\Delta y_{l}^{4}=\Delta\gamma_{1}^{\left(  kl\right)  }=0,\\
a_{21}^{\left(  kl\right)  }\Delta y_{l}=\Delta\gamma_{2}^{\left(  kl\right)
}-a_{22}^{\left(  kl\right)  }\Delta y_{l}^{2}+a_{23}^{\left(  kl\right)
}\Delta y_{l}^{3}+a_{24}^{\left(  kl\right)  }\Delta y_{l}^{4},\\
a_{31}^{\left(  kl\right)  }\Delta y_{l}=\Delta\gamma_{3}^{\left(  kl\right)
}-a_{32}^{\left(  kl\right)  }\Delta y_{l}^{2}+a_{33}^{\left(  kl\right)
}\Delta y_{l}^{3}+a_{34}^{\left(  kl\right)  }\Delta y_{l}^{4},\\
a_{41}^{\left(  kl\right)  }\Delta y_{l}=\Delta\gamma_{4}^{\left(  kl\right)
}-a_{42}^{\left(  kl\right)  }\Delta y_{l}^{2}+a_{43}^{\left(  kl\right)
}\Delta y_{l}^{3}+a_{44}^{\left(  kl\right)  }\Delta y_{l}^{4},\\
a_{11}^{\left(  kl\right)  }\Delta x_{k}+a_{21}^{\left(  kl\right)  }\Delta
x_{k}^{2}+a_{31}^{\left(  kl\right)  }\Delta x_{k}^{3}+a_{41}^{\left(
kl\right)  }\Delta x_{k}^{4}=\Delta\beta_{1}^{\left(  kl\right)  }=0,\\
a_{12}^{\left(  kl\right)  }\Delta x_{k}=\Delta\beta_{2}^{\left(  kl\right)
}-a_{22}^{\left(  kl\right)  }\Delta x_{k}^{2}-a_{32}^{\left(  kl\right)
}\Delta x_{k}^{3}-a_{42}^{\left(  kl\right)  }\Delta x_{k}^{4},\\
a_{13}^{\left(  kl\right)  }\Delta x_{k}=\Delta\beta_{3}^{\left(  kl\right)
}-a_{23}^{\left(  kl\right)  }\Delta x_{k}^{2}-a_{33}^{\left(  kl\right)
}\Delta x_{k}^{3}-a_{43}^{\left(  kl\right)  }\Delta x_{k}^{4}%
\end{array}
\right.  \label{sys7}%
\end{equation}
with the solution%
\begin{equation}
\left\{
\begin{array}
[c]{l}%
a_{13}^{\left(  kl\right)  }=\Delta\beta_{3}^{\left(  kl\right)  }\Delta
x_{k}^{-1}-a_{23}^{\left(  kl\right)  }\Delta x_{k}-a_{33}^{\left(  kl\right)
}\Delta x_{k}^{2}-a_{43}^{\left(  kl\right)  }\Delta x_{k}^{3},\\
a_{12}^{\left(  kl\right)  }=\Delta\beta_{2}^{\left(  kl\right)  }\Delta
x_{k}^{-1}-a_{22}^{\left(  kl\right)  }\Delta x_{k}-a_{32}^{\left(  kl\right)
}\Delta x_{k}^{2}-a_{42}^{\left(  kl\right)  }\Delta x_{k}^{3},\\
a_{41}^{\left(  kl\right)  }=\Delta\gamma_{4}^{\left(  kl\right)  }\Delta
y_{l}^{-1}-a_{42}^{\left(  kl\right)  }\Delta y_{l}-a_{43}^{\left(  kl\right)
}\Delta y_{l}^{2}-a_{44}^{\left(  kl\right)  }\Delta y_{l}^{3},\\
a_{31}^{\left(  kl\right)  }=\Delta\gamma_{3}^{\left(  kl\right)  }\Delta
y_{l}^{-1}-a_{32}^{\left(  kl\right)  }\Delta y_{l}-a_{33}^{\left(  kl\right)
}\Delta y_{l}^{2}-a_{34}^{\left(  kl\right)  }\Delta y_{l}^{3},\\
a_{21}^{\left(  kl\right)  }=\Delta\gamma_{2}^{\left(  kl\right)  }\Delta
y_{l}^{-1}-a_{22}^{\left(  kl\right)  }\Delta y_{l}-a_{23}^{\left(  kl\right)
}\Delta y_{l}^{2}-a_{24}^{\left(  kl\right)  }\Delta y_{l}^{3},\\
a_{11}^{\left(  kl\right)  }=-a_{21}^{\left(  kl\right)  }\Delta x_{k}%
-a_{31}^{\left(  kl\right)  }\Delta x_{k}^{2}-a_{41}^{\left(  kl\right)
}\Delta x_{k}^{3},\\
a_{14}^{\left(  kl\right)  }=-a_{11}^{\left(  kl\right)  }\Delta y_{l}%
^{-3}-a_{12}^{\left(  kl\right)  }\Delta y_{l}^{-2}-a_{13}^{\left(  kl\right)
}\Delta y_{l}^{-1}.
\end{array}
\right.  \label{koef}%
\end{equation}

\subparagraph{$C^{0}$-spline}

\begin{figure}[tbh]%
\centering
\includegraphics[
height=8.8494cm,
width=11.1101cm
]%
{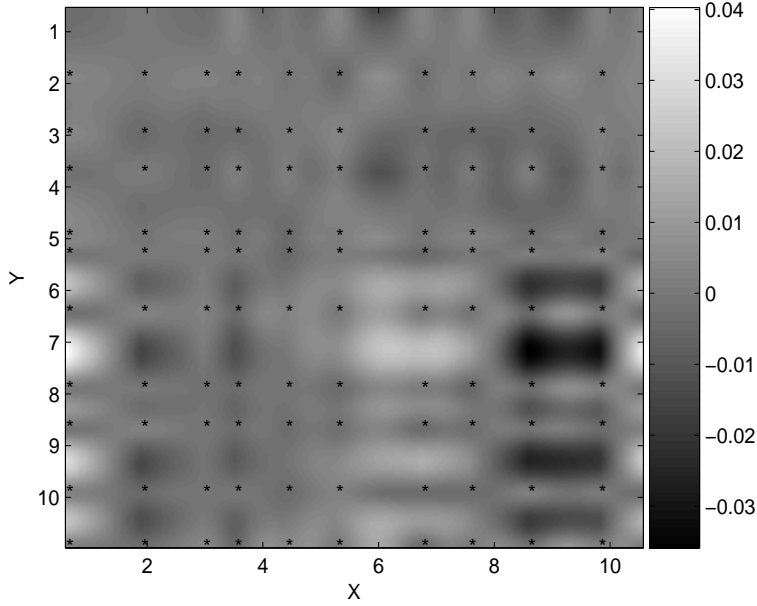}%
\caption{The example of polynomial $C^{0}$-spline with the coefficients
(\ref{koef0}). The points of $\mathfrak{M}$ are marked by stars, $K=L=10$.}%
\label{fig5}%
\end{figure}

Let us set, for example, in (\ref{koef})
\[
a_{ij}^{\left(  kl\right)  }=0,i=\overline{2,4},j=\overline{2,4}.
\]
Then we obtain from (\ref{koef})%
\begin{equation}
\left\{
\begin{array}
[c]{l}%
a_{13}^{\left(  kl\right)  }=\Delta\beta_{3}^{\left(  kl\right)  }\Delta
x_{k},\\
a_{12}^{\left(  kl\right)  }=\Delta\beta_{2}^{\left(  kl\right)  }\Delta
x_{k},\\
a_{41}^{\left(  kl\right)  }=\Delta\gamma_{4}^{\left(  kl\right)  }\Delta
y_{l},\\
a_{31}^{\left(  kl\right)  }=\Delta\gamma_{3}^{\left(  kl\right)  }\Delta
y_{l},\\
a_{21}^{\left(  kl\right)  }=\Delta\gamma_{2}^{\left(  kl\right)  }\Delta
y_{l},\\
a_{11}^{\left(  kl\right)  }=-\left(  \Delta\gamma_{2}^{\left(  kl\right)
}+\Delta\gamma_{3}^{\left(  kl\right)  }\Delta x_{k}+\Delta\gamma_{4}^{\left(
kl\right)  }\Delta x_{k}^{2}\right)  \Delta x_{k}\Delta y_{l},\\
a_{14}^{\left(  kl\right)  }=-a_{11}^{\left(  kl\right)  }\Delta y_{l}%
^{-3}-\Delta\beta_{2}^{\left(  kl\right)  }\Delta x_{k}^{-1}\Delta y_{l}%
^{-2}-\Delta\beta_{3}^{\left(  kl\right)  }\Delta x_{k}^{-1}\Delta y_{l}^{-1}.
\end{array}
\right.  \label{koef0}%
\end{equation}
The spline with such coefficients is continuous but, in general, its first
derivatives have jumps on the boundaries of the cells $\mathbf{P}_{kl}$. The
example of such spline see on Fig. \ref{fig5}%

\subparagraph{$C^{1}$-spline}

Let us try to choose such free parameters that spline (\ref{s2}) becomes
$C^{1}$-function. It is sufficiently to require that the first derivative
$\partial F/\partial x$ is continuous at all vertical lines $x=x_{k}$ and the
$\partial F/\partial y$ is continuous at horizontal lines $y=y_{l}$ of the
grid. This tends to the equations%
\begin{equation}%
\begin{array}
[c]{c}%
a_{1j}^{\left(  kl\right)  }=%
{\displaystyle\sum\limits_{i=1}^{4}}
ia_{ij}^{\left(  k+1\ l\right)  }\Delta x_{k+1}^{i-1},\ a_{i1}^{\left(
kl\right)  }=%
{\displaystyle\sum\limits_{j=1}^{4}}
ja_{ij}^{\left(  k\ l+1\right)  }\Delta y_{l+1}^{j-1},\\
j=\overline{1,4},\ i=\overline{2,4},k=\overline{1,K-1},~l=\overline{1,L-1}.
\end{array}
\label{a1}%
\end{equation}
We see that the second row and the second column of the matrix are completely
determined by the coefficients of the matrix from contiguous cells (that lie
to right and above from a given cell). On the other hand due to the continuity
of the spline and the equations (\ref{koef}) these elements can be found by
nine coefficients of the matrix%
\begin{equation}
\left\Vert
\begin{array}
[c]{ccc}%
a_{22}^{\left(  kl\right)  } & a_{23}^{\left(  kl\right)  } & a_{24}^{\left(
kl\right)  }\\
a_{32}^{\left(  kl\right)  } & a_{33}^{\left(  kl\right)  } & a_{34}^{\left(
kl\right)  }\\
a_{42}^{\left(  kl\right)  } & a_{43}^{\left(  kl\right)  } & a_{44}^{\left(
kl\right)  }%
\end{array}
\right\Vert \label{a22}%
\end{equation}
From the first and the fifth equations of the (\ref{sys7}) and from the
(\ref{a1}) we receive%
\[%
\begin{array}
[c]{c}%
{\displaystyle\sum\limits_{j=1}^{4}}
a_{1j}^{\left(  kl\right)  }\Delta y_{l}^{j-1}=%
{\displaystyle\sum\limits_{j=1}^{4}}
{\displaystyle\sum\limits_{i=1}^{4}}
ia_{ij}^{\left(  k+1\ l\right)  }\Delta x_{k+1}^{i-1}\Delta y_{l}^{j-1}=0,\\%
{\displaystyle\sum\limits_{i=1}^{4}}
a_{i1}^{\left(  kl\right)  }\Delta x_{k}^{i-1}=%
{\displaystyle\sum\limits_{i=1}^{4}}
{\displaystyle\sum\limits_{j=1}^{4}}
ja_{ij}^{\left(  k\ l+1\right)  }\Delta y_{l+1}^{j-1}\Delta x_{k}^{i-1}=0.
\end{array}
\]
So we obtain for $a_{ij}^{\left(  kl\right)  }$%

\begin{equation}
\left\{
\begin{array}
[c]{l}%
{\displaystyle\sum\limits_{i=1}^{4}}
{\displaystyle\sum\limits_{j=1}^{4}}
ja_{ij}^{\left(  kl\right)  }\Delta y_{l}^{j-1}\Delta x_{k}^{i-1}=0,\\%
{\displaystyle\sum\limits_{j=1}^{4}}
{\displaystyle\sum\limits_{i=1}^{4}}
ia_{ij}^{\left(  kl\right)  }\Delta x_{k}^{i-1}\Delta y_{l}^{j-1}=0.
\end{array}
\right.  \label{pkl}%
\end{equation}
Five additional equations we obtain from the system (\ref{sys7}),%
\begin{equation}
\left\{
\begin{array}
[c]{l}%
a_{22}^{\left(  kl\right)  }+a_{32}^{\left(  kl\right)  }\Delta x_{k}%
+a_{42}^{\left(  kl\right)  }\Delta x_{k}^{2}=\Delta\beta_{2}^{\left(
kl\right)  }\Delta x_{k}^{-2}-a_{12}^{\left(  kl\right)  }\Delta x_{k}^{-1},\\
a_{23}^{\left(  kl\right)  }+a_{33}^{\left(  kl\right)  }\Delta x_{k}%
+a_{43}^{\left(  kl\right)  }\Delta x_{k}^{2}=\Delta\beta_{3}^{\left(
kl\right)  }\Delta x_{k}^{-2}-a_{13}^{\left(  kl\right)  }\Delta x_{k}^{-1},\\
a_{22}^{\left(  kl\right)  }+a_{23}^{\left(  kl\right)  }\Delta y_{l}%
+a_{24}^{\left(  kl\right)  }\Delta y_{l}^{2}=\Delta\gamma_{2}^{\left(
kl\right)  }\Delta y_{l}^{-2}-a_{21}^{\left(  kl\right)  }\Delta y_{l}^{-1},\\
a_{32}^{\left(  kl\right)  }+a_{33}^{\left(  kl\right)  }\Delta y_{l}%
+a_{34}^{\left(  kl\right)  }\Delta y_{l}^{2}=\Delta\gamma_{3}^{\left(
kl\right)  }\Delta y_{l}^{-2}-a_{31}^{\left(  kl\right)  }\Delta y_{l}^{-1},\\
a_{42}^{\left(  kl\right)  }+a_{43}^{\left(  kl\right)  }\Delta y_{l}%
+a_{44}^{\left(  kl\right)  }\Delta y_{l}^{2}=\Delta\gamma_{4}^{\left(
kl\right)  }\Delta y_{l}^{-2}-a_{41}^{\left(  kl\right)  }\Delta y_{l}^{-1}.
\end{array}
\right.  \label{s5}%
\end{equation}
From systems (\ref{s5}) and (\ref{pkl}) we receive the system of seven
equations for nine coefficients (\ref{a22})%
\begin{equation}
\left\{
\begin{array}
[c]{l}%
a_{22}^{\left(  kl\right)  }+a_{32}^{\left(  kl\right)  }\Delta x_{k}%
+a_{42}^{\left(  kl\right)  }\Delta x_{k}^{2}=\Delta\beta_{2}^{\left(
kl\right)  }\Delta x_{k}^{-2}-a_{12}^{\left(  kl\right)  }\Delta x_{k}^{-1},\\
a_{23}^{\left(  kl\right)  }+a_{33}^{\left(  kl\right)  }\Delta x_{k}%
+a_{43}^{\left(  kl\right)  }\Delta x_{k}^{2}=\Delta\beta_{3}^{\left(
kl\right)  }\Delta x_{k}^{-2}-a_{13}^{\left(  kl\right)  }\Delta x_{k}^{-1},\\
a_{22}^{\left(  kl\right)  }+a_{23}^{\left(  kl\right)  }\Delta y_{l}%
+a_{24}^{\left(  kl\right)  }\Delta y_{l}^{2}=\Delta\gamma_{2}^{\left(
kl\right)  }\Delta y_{l}^{-2}-a_{21}^{\left(  kl\right)  }\Delta y_{l}^{-1},\\
a_{32}^{\left(  kl\right)  }+a_{33}^{\left(  kl\right)  }\Delta y_{l}%
+a_{34}^{\left(  kl\right)  }\Delta y_{l}^{2}=\Delta\gamma_{3}^{\left(
kl\right)  }\Delta y_{l}^{-2}-a_{31}^{\left(  kl\right)  }\Delta y_{l}^{-1},\\
a_{42}^{\left(  kl\right)  }+a_{43}^{\left(  kl\right)  }\Delta y_{l}%
+a_{44}^{\left(  kl\right)  }\Delta y_{l}^{2}=\Delta\gamma_{4}^{\left(
kl\right)  }\Delta y_{l}^{-2}-a_{41}^{\left(  kl\right)  }\Delta y_{l}^{-1},\\%
{\displaystyle\sum\limits_{i=2}^{4}}
{\displaystyle\sum\limits_{j=2}^{4}}
ja_{ij}^{\left(  kl\right)  }\Delta y_{l}^{j-1}\Delta x_{k}^{i-1}=-%
{\displaystyle\sum\limits_{j=1}^{4}}
ja_{1j}^{\left(  kl\right)  }\Delta y_{l}^{j-1}-%
{\displaystyle\sum\limits_{i=2}^{4}}
a_{i1}^{\left(  kl\right)  }\Delta x_{k}^{i-1},\\%
{\displaystyle\sum\limits_{j=2}^{4}}
{\displaystyle\sum\limits_{i=2}^{4}}
ia_{ij}^{\left(  kl\right)  }\Delta y_{l}^{j-1}\Delta x_{k}^{i-1}=-%
{\displaystyle\sum\limits_{j=1}^{4}}
a_{1j}^{\left(  kl\right)  }\Delta y_{l}^{j-1}-%
{\displaystyle\sum\limits_{i=2}^{4}}
ia_{i1}^{\left(  kl\right)  }\Delta x_{k}^{i-1}.
\end{array}
\right.  \label{sys7a}%
\end{equation}
However the matrix of this system%
\[
\left\Vert
\begin{array}
[c]{ccccccc}%
1 & 0 & 1 & 0 & 0 & 2\Delta y_{l}\Delta x_{k} & 2\Delta x_{k}\Delta y_{l}\\
\Delta x_{k} & 0 & 0 & 1 & 0 & 2\Delta y_{l}\Delta x_{k}^{2} & 3\Delta
x_{k}^{2}\Delta y_{l}\\
\Delta x_{k}^{2} & 0 & 0 & 0 & 1 & 2\Delta y_{l}\Delta x_{k}^{3} & 4\Delta
x_{k}^{3}\Delta y_{l}\\
0 & 1 & \Delta y_{l} & 0 & 0 & 3\Delta y_{l}^{2}\Delta x_{k} & 2\Delta
x_{k}\Delta y_{l}^{2}\\
0 & \Delta x_{k} & 0 & \Delta y_{l} & 0 & 3\Delta y_{l}^{2}\Delta x_{k}^{2} &
3\Delta x_{k}^{2}\Delta y_{l}^{2}\\
0 & \Delta x_{k}^{2} & 0 & 0 & \Delta y_{l} & 3\Delta y_{l}^{2}\Delta
x_{k}^{3} & 4\Delta x_{k}^{3}\Delta y_{l}^{2}\\
0 & 0 & \Delta y_{l}^{2} & 0 & 0 & 4\Delta y_{l}^{3}\Delta x_{k} & 2\Delta
x_{k}\Delta y_{l}^{3}\\
0 & 0 & 0 & \Delta y_{l}^{2} & 0 & 4\Delta y_{l}^{3}\Delta x_{k}^{2} & 3\Delta
x_{k}^{2}\Delta y_{l}^{3}\\
0 & 0 & 0 & 0 & \Delta y_{l}^{2} & 4\Delta y_{l}^{3}\Delta x_{k}^{3} & 4\Delta
x_{k}^{3}\Delta y_{l}^{3}%
\end{array}
\right\Vert
\]
$\allowbreak\allowbreak$has the rank $5$. This means that system (\ref{sys7a})
has no solution. Hence two-dimensional problem \textbf{B }has no solution in
the class of $C^{1}$-splines.

Besides the problem of smoothness the next problem is the appearance of false
stationary points that do not belong to $\mathfrak{M}$. The above
considerations show that class of polynomial splines is not convenient for
solutions of the problem \textbf{B,C}.

\subsection{Trigonometric splines}

\subsubsection{One-dimensional problem}

\label{spline1t}

Consider the following problem \textbf{B}: \emph{for given sets }$\left(
x_{i},z_{i}\right)  $\emph{, }$a=x_{0}<x_{1}<\ldots<x_{N}=b$\emph{ build a
trigonometric spline }$F\left(  \cdot\right)  \in C^{1}\left[  a,b\right]
$\emph{\ such that}%
\begin{equation}
\left\{
\begin{array}
[c]{l}%
F\left(  x\right)  =t_{i}\left(  x\right)  ,\ x\in\left[  x_{i-1}%
,x_{i}\right]  ,\ i=\overline{1,N},\\
t_{i}\left(  x\right)  \equiv a_{i}+b_{i}\cos\left(  c_{i}x+d_{i}\right)  ,\\
F\left(  x_{i}\right)  =z_{i},\\
F^{\prime}\left(  x_{i}\right)  =0.
\end{array}
\right.  \label{trig1}%
\end{equation}
Conditions (\ref{trig1}) do not identify the spline. For example, we can take
the following coefficients
\[
\left\{
\begin{array}
[c]{l}%
a_{i}=\left(  z_{i-1}+z_{i}\right)  /2,\\
b_{i}=\left(  z_{i-1}-z_{i}\right)  /2,\\
c_{i}=\pi/\bigtriangleup x_{i},\\
d_{i}=-\pi x_{i-1}/\bigtriangleup x_{i},\\
\bigtriangleup x_{i}\equiv x_{i}-x_{i-1}.
\end{array}
\right.
\]
The corresponding spline has the form%
\[
t_{i}\left(  x\right)  =\frac{z_{i-1}+z_{i}}{2}+\frac{z_{i-1}-z_{i}}{2}\cos
\pi\frac{x-x_{i-1}}{\bigtriangleup x_{i}}.
\]
For such coefficients false stationary points are absent; this means that
really the spline solves the problem \textbf{C}.

In general, for arbitrary values $z_{i}$ the spline has jumps at points
$x_{i}$. If to forbid the jumps,%
\[%
\begin{array}
[c]{c}%
t_{i-1}^{\prime\prime}\left(  x_{i-1}+0\right)  =-\pi^{2}\left(
z_{i-2}-z_{i-1}\right)  \bigtriangleup x_{i-1}^{-2}/2=\\
=\pi^{2}\left(  z_{i-1}-z_{i}\right)  \bigtriangleup x_{i}^{-2}/2=t_{i}%
^{\prime\prime}\left(  x_{i-1}-0\right)  ,
\end{array}
\]
then we receive the relations%
\[
\left(  z_{i}-z_{i-1}\right)  /\bigtriangleup x_{i}^{2}=\mu\left(  -1\right)
^{i},\ i=\overline{1,n},
\]
where $\mu$ is an arbitrary value. Such kind $Ñ ^{2}$-splines form
two-parametrical set as for fixed value $z_{0}$ other values $z_{i}$ can be
calculated by the formula%
\begin{equation}
z_{i}=z_{0}+\mu%
{\displaystyle\sum\limits_{l=1}^{i}}
\left(  -1\right)  ^{l}\bigtriangleup x_{l}^{2}\label{trig12}%
\end{equation}
which contains two free parameters $z_{0}$ and $\mu.$ For fixed set $\left\{
x_{i}\right\}  $ the parameter $\mu$ evaluates mean deviation of the spline
from the initial value $z_{0}$:%
\[%
{\displaystyle\sum\limits_{i=1}^{n}}
\left(  z_{i}-z_{0}\right)  /n=\mu%
{\displaystyle\sum\limits_{l=1}^{n}}
\left(  -1\right)  ^{l}\left(  1-l/n\right)  \bigtriangleup x_{l}^{2}.
\]

\begin{figure}[tbh]%
\centering
\includegraphics[
height=7.7668cm,
width=9.42cm
]%
{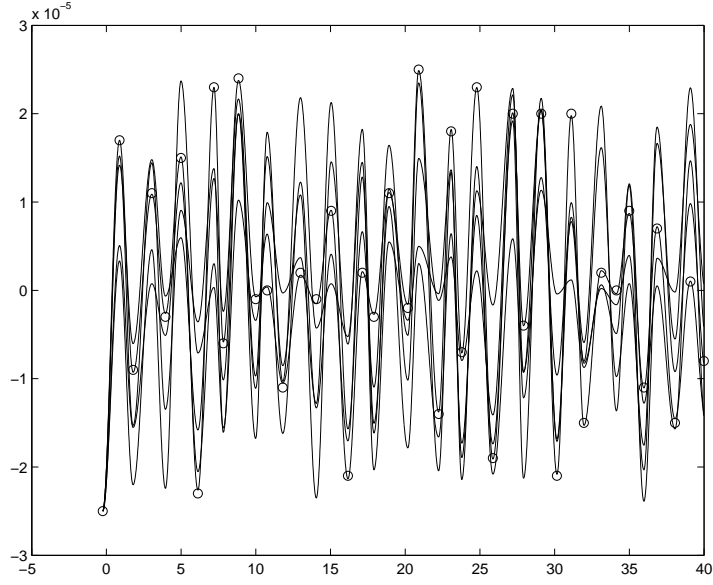}%
\caption{A beam of smooth trigonometric splines. The values $\left\{
z_{i}\right\}  $ were  taken randomly from a given interval.}%
\label{figtrig11}%
\end{figure}

Examples of trigonometric splines (\ref{trig1}) are presented on Fig.
\ref{figtrig11}%

\subsubsection{Two-dimensional problem}

\label{spline2t}

Consider the following problem \textbf{B}: \emph{for given sets }$\left(
x_{k},y_{l},z_{kl}\right)  $\emph{, }$a=x_{0}<x_{1}<\ldots<x_{K}=b$,
$c=y_{0}<y_{1}<\ldots<y_{L}=d$\emph{ build a spline }$F\left(  \cdot\right)
\in C^{1}\left(  \left[  a,b\right]  \times\left[  c,d\right]  \right)
$\emph{\ such that}%

\begin{equation}
\left\{
\begin{array}
[c]{l}%
F\left(  x,y\right)  =t_{kl}\left(  x,y\right)  ,\ \left(  x,y\right)
\in\mathbf{P}_{kl},\ k=\overline{1,K},~l=\overline{1,L},\\
t_{kl}\left(  x,y\right)  =\left[  a_{kl}+b_{kl}\cos\left(  c_{kl}%
x+d_{kl}\right)  \right]  \left[  A_{kl}+B_{kl}\cos\left(  C_{kl}%
y+D_{kl}\right)  \right]  ,\\
F\left(  x_{k},y_{l}\right)  =z_{kl},\\
\triangledown F\left(  x_{k},y_{l}\right)  =0.
\end{array}
\right.  \label{trig2}%
\end{equation}

In analogues to one-dimensional case the conditions (\ref{trig2}) do not
identify the spline. For example, we can take%
\begin{equation}
\left\{
\begin{array}
[c]{l}%
c_{kl}=\pi/\Delta x_{k},\\
d_{kl}=-\pi x_{k-1}/\Delta x_{k},\\
C_{kl}=\pi/\Delta y_{l},\\
D_{kl}=-\pi y_{l-1}/\Delta y_{l},
\end{array}
\right.  \label{coeft1}%
\end{equation}
where $\Delta x_{k}\equiv x_{k}-x_{k-1}$, $\Delta y_{l}\equiv y_{l}-y_{l-1}$.
In this case $\triangledown F\left(  x_{k},y_{l}\right)  =0$ at the vertex of
$\mathbf{P}_{kl}$ (and only at these points; so, actually the spline solves
the problem \textbf{C}). For other coefficients of the spline we have the
equations%
\[
\left\{
\begin{array}
[c]{l}%
t_{kl}\left(  x_{k-1},y_{l-1}\right)  =\left[  a_{kl}+b_{kl}\right]  \left[
A_{kl}+B_{kl}\right]  =z_{k-1l-1},\\
t_{kl}\left(  x_{k},y_{l-1}\right)  =\left[  a_{kl}-b_{kl}\right]  \left[
A_{kl}+B_{kl}\right]  =z_{kl-1},\\
t_{kl}\left(  x_{k-1},y_{l}\right)  =\left[  a_{kl}+b_{kl}\right]  \left[
A_{kl}-B_{kl}\right]  =z_{k-1l},\\
t_{kl}\left(  x_{k},y_{l}\right)  =\left[  a_{kl}-b_{kl}\right]  \left[
A_{kl}-B_{kl}\right]  =z_{kl}%
\end{array}
\right.
\]
which give%
\begin{equation}
\left\{
\begin{array}
[c]{l}%
a_{kl}A_{kl}=\left(  z_{k-1l-1}+z_{k-1l}+z_{kl-1}+z_{kl}\right)  /4,\\
b_{kl}A_{kl}=\left(  z_{k-1l-1}+z_{k-1l}-z_{kl-1}-z_{kl}\right)  /4,\\
a_{kl}B_{kl}=\left(  z_{k-1l-1}+z_{kl-1}-z_{k-1l}-z_{kl}\right)  /4,\\
b_{kl}B_{kl}=\left(  z_{k-1l-1}+z_{kl}-z_{k-1l}-z_{kl-1}\right)  /4,
\end{array}
\right.  \label{coeft2}%
\end{equation}
that is sufficient for the spline determination. As coefficients of
$t_{kl}\left(  x,y\right)  $ depend only on the spline values at vertexes of
$\mathbf{P}_{kl}$, then we must restrict these values in order to provide the
smoothness of the spline.

It is easy to prove that we need to provide the smoothness of the spline on
the boundaries of $\mathbf{P}_{kl}$ only. Take for example the common boundary
of $\mathbf{P}_{kl}$ and $\mathbf{P}_{kl-1}$, the segment $\left[
x_{k-1},x_{k}\right]  \times\left\{  y_{l-1}\right\}  $. The spline is
continuous on this segment if and only if for all $x\in\left[  x_{k-1}%
,x_{k}\right]  $%
\[%
\begin{array}
[c]{c}%
\left[  a_{kl}+b_{kl}\cos\pi\left(  x-x_{k-1}\right)  /\Delta x_{k}\right]
\left[  A_{kl}+B_{kl}\right]  =\\
=\left[  a_{kl-1}+b_{kl-1}\cos\pi\left(  x-x_{k-1}\right)  /\Delta
x_{k}\right]  \left[  A_{kl-1}-B_{kl-1}\right]  .
\end{array}
\]
So we have the equations%
\begin{equation}
\left\{
\begin{array}
[c]{l}%
a_{kl-1}\left(  A_{kl-1}-B_{kl-1}\right)  =a_{kl}\left(  A_{kl}+B_{kl}\right)
,\\
b_{kl-1}\left(  A_{kl-1}-B_{kl-1}\right)  =b_{kl}\left(  A_{kl}+B_{kl}\right)
.
\end{array}
\right.  \label{contt2}%
\end{equation}

Substitute formulae (\ref{coeft2}) to these equations we receive the
identities for $z_{kl}$. Also we get the identities when check the continuity
of the spline on the vertical boundary $\left\{  x_{k-1}\right\}
\times\left[  y_{l-1},y_{l}\right]  $. Hence the spline is continuous for all
values $\left\{  z_{kl}\right\}  $. Actually the spline belongs even to the
class $C^{1}\left(  \left[  a,b\right]  \times\left[  c,d\right]  \right)  $.
To prove it we must prove the continuity of the gradient on the boundaries of
$\mathbf{P}_{kl}$. For example, using (\ref{contt2}), we have on the segment
$\left[  x_{k-1},x_{k}\right]  \times\left\{  y_{l-1}\right\}  $
\[%
\begin{array}
[c]{l}%
F_{y}^{\prime}\left(  x,y_{l-1}\pm0\right)  \equiv0,\\
F_{x}^{\prime}\left(  x,y_{l-1}\pm0\right)  \equiv\pi b_{kl}\sin\left(
\pi\left(  x-x_{k-1}\right)  /\Delta x_{k}\right)  \left[  A_{kl}%
+B_{kl}\right]  /\Delta x_{k}\equiv\\
\equiv\pi b_{kl-1}\sin\left(  \pi\left(  x-x_{k-1}\right)  /\Delta
x_{k}\right)  \left[  A_{kl-1}-B_{kl-1}\right]  /\Delta x_{k}.
\end{array}
\]

If we want to construct the more smooth spline then we must to impose its
values $\left\{  z_{kl}\right\}  $. For example the condition of continuity of
the second derivative $F_{y^{2}}^{\prime\prime}$ on the boundary $\left[
x_{k-1},x_{k}\right]  \times\left\{  y_{l-1}\right\}  $ tends to the equations%
\[%
\begin{array}
[c]{c}%
\left[  a_{kl}+b_{kl}\cos\pi\left(  x-x_{k-1}\right)  /\Delta x_{k}\right]
B_{kl}C_{kl}^{2}=\\
=-\left[  a_{kl-1}+b_{kl-1}\cos\pi\left(  x-x_{k-1}\right)  /\Delta
x_{k}\right]  B_{kl-1}C_{kl-1}^{2},
\end{array}
\]
hence,%
\[
a_{kl}B_{kl}C_{kl}^{2}=-a_{kl-1}B_{kl-1}C_{kl-1}^{2},\ b_{kl}B_{kl}C_{kl}%
^{2}=-b_{kl-1}B_{kl-1}C_{kl-1}^{2}.
\]
If we put (\ref{coeft2}) to these equalities, then we get%
\begin{equation}%
\begin{array}
[c]{c}%
\left(  z_{k-1l-1}+z_{kl-1}-z_{k-1l}-z_{kl}\right)  \beta_{l}=\left(
z_{k-1l-1}+z_{kl-1}-z_{k-1l-2}-z_{kl-2}\right)  ,\\
\left(  z_{k-1l-1}+z_{kl}-z_{k-1l}-z_{kl-1}\right)  \beta_{l}=\left(
z_{k-1l-1}+z_{kl-2}-z_{k-1l-2}-z_{kl-1}\right)  ,
\end{array}
\label{zyy}%
\end{equation}
where $\beta_{l}\equiv\left(  \Delta y_{l-1}/\Delta y_{l}\right)  ^{2}$. In
analogues the condition of continuity of the second derivative $F_{x^{2}%
}^{\prime\prime}$ on the boundary$\left\{  x_{k-1}\right\}  \times\left[
y_{l-1},y_{l}\right]  $ tends to the equations%
\[%
\begin{array}
[c]{c}%
\left[  A_{kl}+B_{kl}\cos\pi\left(  y-y_{l-1}\right)  /\Delta y_{l}\right]
b_{kl}c_{kl}^{2}=\\
=-\left[  A_{k-1l}+B_{k-1l}\cos\pi\left(  y-y_{l-1}\right)  /\Delta
y_{l}\right]  b_{k-1l}c_{k-1l}^{2},
\end{array}
\]
hence,
\[
A_{kl}b_{kl}c_{kl}^{2}=-A_{k-1l}b_{k-1l}c_{k-1l}^{2},\ B_{kl}b_{kl}c_{kl}%
^{2}=-B_{k-1l}b_{k-1l}c_{k-1l}^{2},
\]
and
\begin{equation}%
\begin{array}
[c]{c}%
\left(  z_{k-1l-1}+z_{k-1l}-z_{kl-1}-z_{kl}\right)  \delta_{k}=z_{k-1l-1}%
+z_{k-1l}-z_{k-2l-1}-z_{k-2l},\\
\left(  z_{k-1l-1}+z_{kl}-z_{k-1l}-z_{kl-1}\right)  \delta_{k}=z_{k-2l}%
+z_{k-1l-1}-z_{k-2l-1}-z_{k-1l},
\end{array}
\label{zxx}%
\end{equation}
where $\delta_{k}\equiv\left(  \Delta x_{k-1}/\Delta x_{k}\right)  ^{2}$. It
is easy to prove that other derivatives of the second order are also
continuous in this case.

From (\ref{zyy}) and (\ref{zxx}) we have%
\[
\left\{
\begin{array}
[c]{l}%
z_{kl-2}=\beta_{l}z_{kl}+z_{kl-1}\left(  1-\beta_{l}\right)  ,\\
z_{k-1l-2}=\beta_{l}z_{k-1l}+z_{k-1l-1}\left(  1-\beta_{l}\right)  ,\\
z_{k-2l}=\delta_{k}z_{kl}+z_{k-1l}\left(  1-\delta_{k}\right)  ,\\
z_{k-2l-1}=\delta_{k}z_{kl-1}+z_{k-1l-1}\left(  1-\delta_{k}\right)  .
\end{array}
\right.
\]
Differently these formulas can be written in the form%
\[%
\begin{array}
[c]{l}%
\left(  z_{kl-1}-z_{kl-2}\right)  /\Delta y_{l-1}^{2}=-\left(  z_{kl}%
-z_{kl-1}\right)  /\Delta y_{l}^{2},\\
\left(  z_{k-1l}-z_{k-2l}\right)  /\Delta x_{k-1}^{2}=-\left(  z_{kl}%
-z_{k-1l}\right)  /\Delta x_{k}^{2},
\end{array}
\]
and this means that%
\[%
\begin{array}
[c]{c}%
\left(  z_{kl}-z_{kl-1}\right)  /\Delta y_{l}^{2}=\left(  -1\right)  ^{l}%
\mu_{k},\ l=\overline{1,L},\ k=\overline{0,K}\\
\left(  z_{kl}-z_{k-1l}\right)  /\Delta x_{k}^{2}=\left(  -1\right)  ^{k}%
\nu_{l},\ k=\overline{1,K},\ l=\overline{0,L}%
\end{array}
\]
for some sets of values $\left\{  \mu_{k}\right\}  ,$ $\left\{  \nu
_{l}\right\}  $. So,%
\[
z_{kl}=z_{k0}+\mu_{k}%
{\displaystyle\sum\limits_{j=1}^{l}}
\left(  -1\right)  ^{j}\Delta y_{j}^{2}=z_{0l}+\nu_{l}%
{\displaystyle\sum\limits_{i=1}^{k}}
\left(  -1\right)  ^{i}\Delta x_{i}^{2},
\]%
\[
z_{kl}-z_{k-1l}=z_{k0}-z_{k-10}+\left(  \mu_{k}-\mu_{k-1}\right)
{\displaystyle\sum\limits_{j=1}^{l}}
\left(  -1\right)  ^{j}\Delta y_{j}^{2}=\nu_{l}\left(  -1\right)  ^{k}\Delta
x_{k}^{2},
\]%
\[
\left(  -1\right)  ^{k}\frac{\left(  \mu_{k}-\mu_{k-1}\right)  }{\Delta
x_{k}^{2}}%
{\displaystyle\sum\limits_{j=1}^{l}}
\left(  -1\right)  ^{j}\Delta y_{j}^{2}=\nu_{l}-\nu_{0}.
\]
For fulfilment of the last equalities for all $k,l$ one can take for any value
$\lambda$%
\[
\mu_{k}=\mu_{0}+\lambda%
{\displaystyle\sum\limits_{j=1}^{k}}
\left(  -1\right)  ^{j}\Delta x_{j}^{2}.
\]
In this case%
\[
\nu_{l}=\nu_{0}+\lambda%
{\displaystyle\sum\limits_{j=1}^{l}}
\left(  -1\right)  ^{j}\Delta y_{j}^{2},
\]%
\begin{equation}%
\begin{array}
[c]{c}%
z_{kl}=z_{00}+\nu_{0}%
{\displaystyle\sum\limits_{i=1}^{k}}
\left(  -1\right)  ^{i}\Delta x_{i}^{2}+\mu_{0}%
{\displaystyle\sum\limits_{j=1}^{l}}
\left(  -1\right)  ^{j}\Delta y_{j}^{2}+\\
+\lambda%
{\displaystyle\sum\limits_{i=1}^{k}}
\left(  -1\right)  ^{i}\Delta x_{i}^{2}%
{\displaystyle\sum\limits_{j=1}^{l}}
\left(  -1\right)  ^{j}\Delta y_{j}^{2}.
\end{array}
\label{trig22}%
\end{equation}
We see that this family of $Ñ ^{2}$-splines depends only on the next four free
parameters $z_{00},\nu_{0},\mu_{0},\lambda.$

\begin{figure}[tbh]%
\centering
\includegraphics[
height=7.7647cm,
width=9.7476cm
]%
{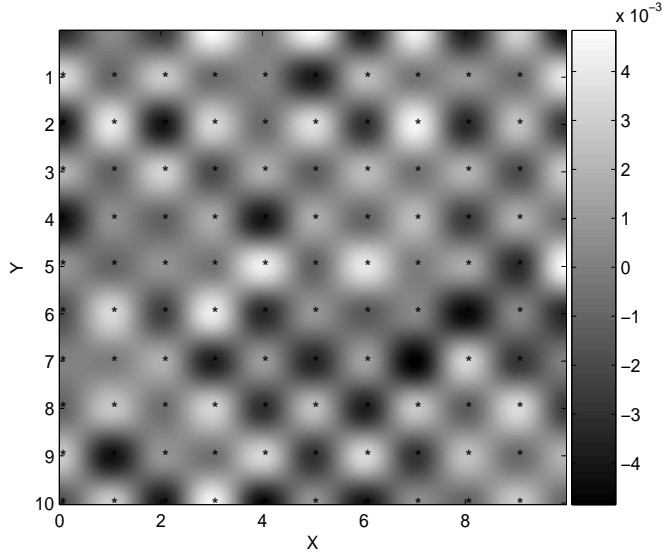}%
\caption{An example of two-dimensional trigonometric spline (\ref{trig2}).
Values $\left\{  z_{ij}\right\}  ,$ $i,j=\overline{1,10}$ were taken randomly
from a fixed segment.}%
\label{figtrig21}%
\end{figure}

\begin{figure}[ptb]%
\centering
\includegraphics[
height=2.8178in,
width=3.6848in
]%
{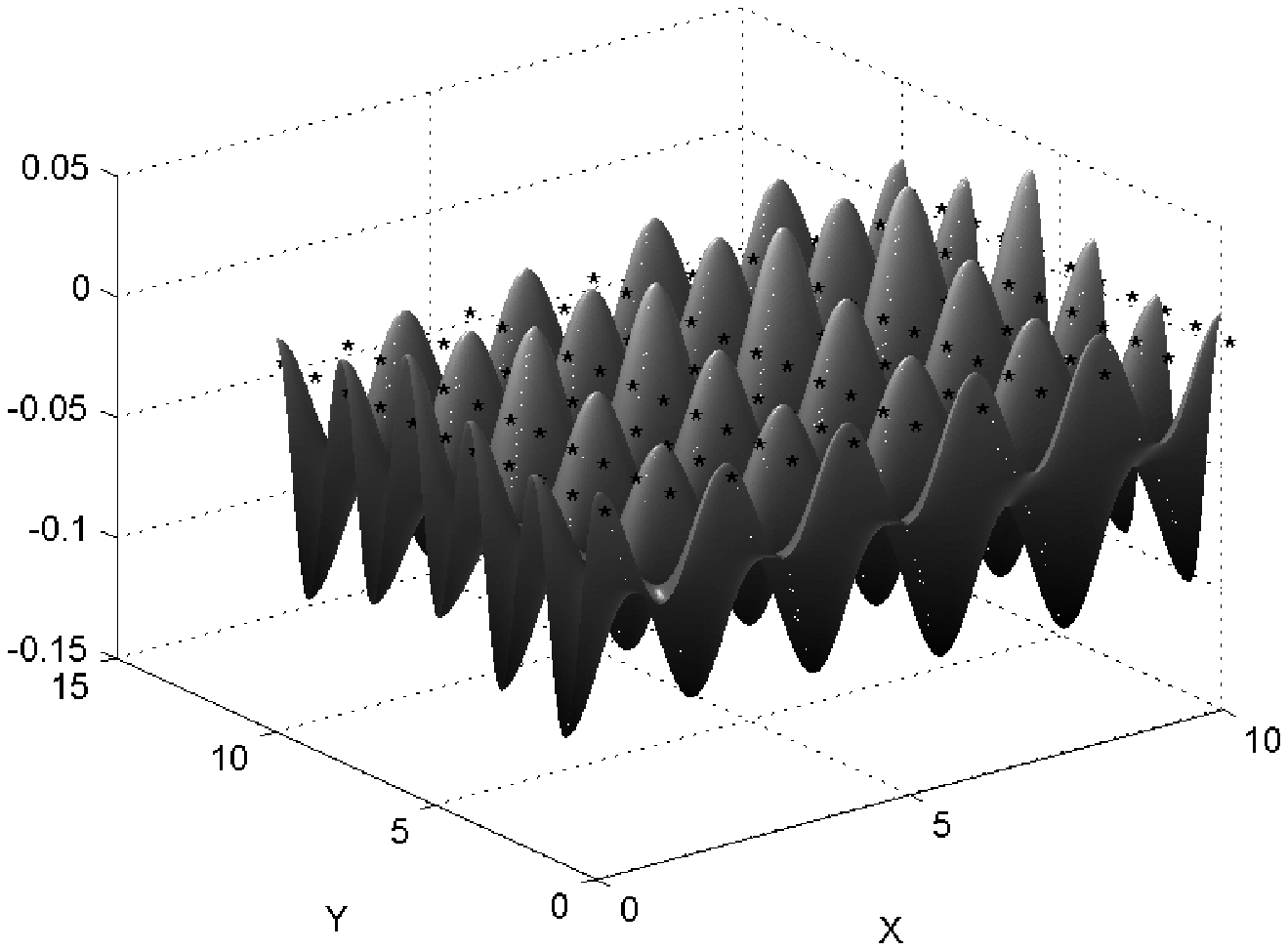}%
\caption{An example of two-dimensional trigonometric $C^{2}$-spline
(\ref{trig2}). Values $\left\{  z_{ij}\right\}  ,$ $i,j=\overline{1,10}$ were
calculated by the formula ( \ref{trig22}). }%
\label{figtrig22}%
\end{figure}

Examples of two-dimensional trigonometric splines see in Figs. \ref{figtrig21}%
,\ref{figtrig22}.%


\section{Superpositions of splines}

\label{antistat}

In the last section we have established that for two-dimensional sets
$\mathfrak{M}$ with matrix structure (\ref{Matr}) solutions of the problems
\textbf{A}-\textbf{C} can be found in the form of splines (polynomial or
trigonometric with different types of smoothness). If $\mathfrak{M}$ does not
have the matrix form any spline has an additional (false) points where its
gradient equals to zero; so it can not solve the problem \textbf{C}. To
overcome this obstacle some constructions may be proposed.

For example, consider the following procedure consisting of the steps

\begin{description}
\item[Step 1.] Complete $\mathfrak{M}$ to a set $\mathfrak{M}_{XY}$ of the
matrix form and build any spline $F\left(  x,y;\mathfrak{M}_{XY}\right)  $
that solves the corresponding problem \textbf{C};

\item[Step 2.] Turn the coordinate system $XY$ for some ("sufficiently small")
angle $\varphi_{1}$ and build an analogous spline in the coordinates
$X_{1}Y_{1}$,%
\[
F\left(  x_{_{1}},y_{_{1}};\mathfrak{M}_{X_{1}Y_{1}}\right)  \equiv F\left(
x\cos\varphi_{1}-y\sin\varphi_{1},x\sin\varphi_{1}+y\cos\varphi_{1}%
;\mathfrak{M}_{X_{1}Y_{1}}\right)  ;
\]

\item[Step 3.] Repeat the above steps several times for different angles
$\varphi_{k}$; then build the function%
\begin{equation}
\tilde{F}\left(  x,y\right)  =%
{\displaystyle\sum\limits_{k\geq1}}
F\left(  x\cos\varphi_{k}-y\sin\varphi_{k},x\sin\varphi_{k}+y\cos\varphi
_{k};\mathfrak{M}_{X_{k}Y_{k}}\right)  .\label{super}%
\end{equation}

\end{description}

As gradients of all summands are zero at points $\mathfrak{M}$, the function
$\tilde{F}\left(  x,y\right)  $ (superposition of splines) solves the problem
\textbf{A}. At the same time the additional points for matrixes $\mathfrak{M}%
_{X_{k}Y_{k}}$ are not the same (they are removed from their initial places);
so, the probability that gradients of the summary function (\ref{super}) at
additional points are not zero is high. Test calculations verify these assumption.

\end{document}